\documentclass{amsart}
\usepackage{graphicx}
\newtheorem{thm}{Theorem}[section]
\newtheorem{cor}[thm]{Corollary}
\newtheorem{lem}[thm]{Lemma}
\newtheorem{prop}[thm]{Proposition}
\theoremstyle{definition}
\newtheorem{defn}[thm]{Definition}
\theoremstyle{remark}

\numberwithin{equation}{section}

\renewcommand{\>}{\rangle}

\newcommand{\cal}[1]{\mathcal{#1}}
\newcommand{\<}{\langle}
\newcommand{\N}{\mathbb{N}}
\begin{document}
\title[]{On finitely generated models of theories with at most countably many
nonisomorphic finitely generated models}%
\author{Abderezak OULD HOUCINE}%
\address{Institut Camille Jordan,
       Universit\'e Claude Bernard Lyon-1,
       B\^atiment Braconnier,
       21 Avenue Claude Bernard,
       69622 Villeurbanne Cedex, France.}%
\email{ould@math.univ-lyon1.fr}%
\maketitle
\begin{abstract}
We  study  finitely generated models of countable theories, having
at most countably many nonisomorphic finitely generated models. We
introduce a notion of rank of finitely generated models  and we
prove, when $T$ has at most countably many nonisomorphic finitely
generated models, that every finitely generated model has an
ordinal rank.  This rank is used to give a property of finitely
generated models analogue to the Hopf property of groups and also
to give a necessary and sufficient condition for  a finitely
generated model to be prime of its complete theory.  We
investigate   some properties of limit groups of equationally
noetherian groups, in respect to their ranks.
\end{abstract}
\section{Introduction}

Throughout this paper, we let $L$ be a fixed countable first order
language; $L$ is arbitrary but is held fixed to simplify notation.
We shall say that a model $\mathcal{M}$ is \emph{$n$-generated} if
it is generated by $n$ elements; that is, if  there exists an
$n$-tuple $\bar a$ in $\mathcal{M}$ such that for every $y \in
\mathcal{M}$, $y=\tau (\bar a)$  for some term $\tau (\bar x)$ of
$L$. It is worth mentioning that this definition is different from
the one used by A. Pillay in \cite{Pillay1, Pillay2}. A theory $T$
is said \emph{$n$-consistent} if it has an $n$-generated model.

The purpose of the paper is to study  some properties of finitely
generated models of countable theories which have at most
countably many nonisomorphic finitely generated models.  Part of
our interest on this case comes from  the fact that the number of
nonisomorphic finitely generated models of a countable theory is
either at most $\aleph_0$ or equals to $2^{\aleph_0}$ (Theorem
\ref{the-number}). This was also motivated by the fact that
several
  theories satisfy the above property. For instance
 the universal (and thus  the complete theory) of a  linear group, over a commutative noetherian
 ring, e.g, a field,  has at most $\aleph_0$
 nonisomorphic finitely generated models. More generally, the universal theory
 of an equationally noetherian group satisfies the same property \cite{ould-equa}.
 We note
 that linear groups over a commutative noetherian
 ring, e.g, a field,  are equationally noetherian and not all
 equationally noetherian groups are linear
 \cite{rem-mya-alge-geo}.

Given a theory $T$, we consider the following type, in the
language $L(\bar c)=L \cup \{c_1, \cdots, c_n\}$, where $c_1,
\cdots, c_n$ are a new constants symbols,
$$
p_n(y)=\{y \neq  \tau (\bar c)~|~\tau(\bar x) \hbox{ a term in }
L\}.
$$

Then a model $(\mathcal M, \bar a)$ of $T$, in the language
$L(\bar c)$, is generated by $\bar a$ if and only $(\mathcal M,
\bar a)$ omits $p_n$. Thus the class of $n$-generated models of
$T$ is the class of models of $T$, in the language $L(\bar c)$,
which omit $p_n$. We adopt this viewpoint,  and for that raison we
need to omit some types in the class of models omitting types
including $p_n$ (Theorem \ref{omi-type} and Theorem
\ref{omi-type2}).

We begin by showing, in the next section,  that the number
$\alpha(T)$ of nonisomorphic finitely generated models of $T$ is
at most $\aleph_0$ or
 equals to $2^{\aleph_0}$. This  result
 is analogue  to the one known on the number of complete types of
 countable theories. Our aim in section 3 is  to define a rank of $n$-generated models of
$T$, and to prove that when $T$ has at most countably many
nonisomorphic $n$-generated models, every
 $n$-generated model of $T$ has an ordinal rank.  As any
 $n$-generated model of $T$ is determined  by its complete
 type relatively to some $n$-generating tuple, the previous result can be seen as an analogue of
 the Morley rank of types of an $\omega$-stable theory $T$. Nevertheless,  the rank that we define is
 different from the Morley rank, since it uses actually the class
 of  $n$-generated models of $T$ and $T$ is not necessarily
 $\omega$-stable.

 In section 4 we turn to use the rank  to give a property of $n$-generated models of
 countable theories, having at most countably many nonisomorphic finitely
 generated models, which can be seen as an analogue to the Hopf property of groups (Theorem \ref{theo-formule-prime}). Recall that a group $G$ is said to be \emph{Hopfian},
 if any surjective morphism from $G$ to $G$ is an isomorphism. It is known that finitely generated linear groups are Hopfian \cite{LyndonSchupp77}.
 More generally, finitely generated equationally noetherian groups are Hopfian
 \cite{ould-equa}. The property obtained for finitely generated
 models is as follows. For every generating $n$-tuple $\bar a$ of $\cal M$,
 there exists a formula $\phi(\bar x)$ such that $\mathcal M
 \models \phi(\bar a)$ and such that for any generating $n$-tuple $\bar b$ of $\cal
 M$, $(\cal M, \bar a) \cong (\mathcal M, \bar b)$ if and only if
 $\mathcal M \models \phi(\bar b)$. In the same section, we use this to get a necessary and sufficient
 condition for a finitely generated model
 to be prime of its complete theory (Theorem \ref{theorem-prime-model}).

 When $T$ is universal, we define an another rank, more easy to use.
 This will be done in section 4.
 Since the universal theory of an equationally noetheiran group
 $H$
 has  at most countably many nonisomorphic finitely generated
 models, one can attribute a rank to $H$-limit groups and study
 them in this context. We give some results in this direction in section 5.

\section{The number of finitely generated
models.}

For a theory $T$  we denote by $\alpha_n(T)$ (resp. $\alpha(T)$)
the number of nonisomorphic $n$-generated (resp. finitely
generated) models of $T$.

\begin{thm} \label{the-number}
For every theory $T$, $\alpha_n(T) \leq \aleph_0$ or
$\alpha_n(T)=2^{\aleph_0}$. Therefore $\alpha(T) \leq \aleph_0$ or
$\alpha(T) = 2^{\aleph_0}$.
\end{thm}

Before proving the theorem, we need some notions. Let $\Phi$ be a
sentence of $L_{\omega_1 \omega}$ and $L_A$ be a countable
fragment of $L_{\omega_1 \omega}$. We define an $L_A$-$n$-type to
be a set $p$ such that
$$
p=\{\varphi(\bar x) \in L_A~|~\mathcal{M}\models \varphi(\bar a)
\},
$$
for some model $\mathcal{M}$ of $\Phi$ and some tuple $\overline a
\in\mathcal{M}^n$. The set of all $L_A$-$n$-types of $\Phi$ is
denoted by $S_n(L_A,\Phi)$. For more details, the reader is
referred to \cite{Marker, Morley}.
\begin{prop}{\rm\cite[Corollary 2.4]{Morley}}\label{coro-morley} For every countable
fragment $L_A$ the set $S_n(L_A,\Phi)$ is either countable or of
power $2^{\aleph_0}$. \qed
\end{prop}

\noindent \textbf{Proof of Theorem \ref{the-number} }

Let

$$\Phi = \bigwedge_{\varphi \in T} \varphi \wedge \exists \bar x \forall
y (\bigvee_{\tau \in Ter} y=\tau(\bar x)),
$$
where $Ter$ is the set of all terms of $L$ with  free variables
among $\bar x$, and the length of $\bar x$ is $n$. Then $\Phi \in
L_{\omega_1\omega}$ and a model of $T$  is  $n$-generated if and
only if it satisfies $\Phi$. Let $L_A$ be the set of all formulas
of $L$. Then $L_A$ is a countable fragment and by Proposition
\ref{coro-morley}, $S_n(L_A,\Phi)$ is countable or of power
$2^{\aleph_0}$.

If $S_n(L_A,\Phi)$ is countable, then $\Phi$ has at most
$\aleph_0$ models; as any complete type of a model of $\Phi$,
respectively to some generating $n$-tuple, is an $L_A$-$n$-type.
Thus $T$ has at most countably many nonisomorphic $n$-generated
models.

If $S_n(L_A,\Phi)$ is of power $2^{\aleph_0}$, then $\Phi$ has
$2^{\aleph_0}$ models (and thus $T$ has $2^{\aleph_0}$
$n$-generated models), as any model of $\Phi$ realizes at most a
countable number of $L_A$-$n$-types. \qed

\bigskip

There are several examples of theories having at most $\aleph_0$
nonisomorphic finitely generated models. For example abelian group
theory, complete theory of finitely generated linear groups, and
more generally the complete theory (or the universal theory) of a
submodel of an $\omega$-stable model. In particular, the universal
theory of non-abelian free groups has $\aleph_0$ nonisomorphic
finitely generated models, as every non-abelian free group is a
subgroup of $GL_n(\mathbb{C})$ and this last group has a finite
Morley rank. This is also a consequence, as noticed in the
introduction,  of the fact that linear groups are equationally
noetherian. We regroup this remarks in the following proposition.

\begin{prop}
Let $T$ be a countable $\omega$-stable theory and $\mathcal M$ a
model of $T$. Then the universal theory of any submodel of
$\mathcal M$ has at most $\aleph_0$ nonisomorphic finitely
generated models.
\end{prop}

\proof Let $\mathcal N \subseteq \mathcal M$ and denote by
$\Gamma$ the universal theory of $\mathcal N$.  Suppose towards a
contradiction  that $\alpha(\Gamma) > \aleph_0$ and thus
$\alpha_n(\Gamma) > \aleph_0$ for some $n \in \mathbb N^*$. Then,
for any $n$-finitely generated model $\mathcal A$ of $\Gamma$,
generated by $\bar a$,  the type
$$ p_{\bar a}(\mathcal A)=\{
\phi(\bar x)~|~\phi \hbox{ is atomic or negatomic such that }
\mathcal A \models \phi(\bar a)\}
$$
is a consistent type with the universal theory of $\mathcal M$. By
compactness, there exists a model of Th$(\cal M)$ which contains a
copy of every $n$-generated model of $\Gamma$. Therefore, as there
exists $\alpha_n(\Gamma) > \aleph_0$ nonisomorphic $n$-generated
models of $\Gamma$, Th$(\cal M)$ has more than $\aleph_0$ types on
$\emptyset$. A contradiction with the $\omega$-stability of
$\mathcal M$. \qed

\smallskip
 One can extracted from  papers of F.Oger \cite{Oger1,Oger3, Oger2} that for every  $n \in \mathbb{N}^*$, there exists
 a complete theory $T$ of groups  such that $\alpha(T)=n$. G.
 Sabbagh asked the following.

\bigskip
 \noindent \textbf{Problem.} Is there a complete theory $T$ of groups such that
 $\alpha(T)=2^{\aleph_0}$ ?

\section{A Rank}

 We shall define  a
\emph{rank} of $n$-generated models of $T$.  Before proceeding, we
need some notions around omitting  types in some classes of
models. For our purpose, a type is a set of sentences in the
language $L(\bar x)$. Let $\mathcal{K}$ be a class of models.
Given a type $q$, we say that $q$ is \emph{supported over
$\mathcal{K}$} if there exists a formula $\phi(\bar x)$ such that:

\smallskip
$(i)$ some model in $\mathcal{K}$ has a tuple satisfying $\phi$,
and

$(ii)$ in every model in $\mathcal{K}$, each tuple satisfying
$\phi$ realizes $q$,

\smallskip
and in that case we say that $\phi$ \emph{supports} $q$ over $\cal
K$. We say that $q$ is \emph{unsupported} over $\cal K$ if it is
not supported over $\cal K$.

\smallskip
Let $T$ be a theory in $L$ and $P$ a set of types. We denote by
$\mathcal{K}(T|P)$ the class of models of $T$ omitting every $p$
in $P$. The next theorem is a slight refinement of the classical
omitting types theorem.

\begin{thm}\label{omi-type}
Let $P$ and $Q$ be a countable sets of types. If
$\mathcal{K}(T|P)$ is not empty and each $q \in Q$ is unsupported
over $\mathcal{K}(T|P)$, then there exists a countable model in
$\mathcal{K}(T|P)$ which omits every $q \in Q$.

\end{thm}

\proof

Let $T'$ be the set of all sentences of $L$ which are true in
every model in $\mathcal{K}(T|P)$. We claim that over $T'$, all
types of $P$ and all types of $Q$ are with no support; that is,
they are unsupported over the class of models of $T'$.

Let $p \in P$ and suppose that $\phi$ supports $p$ over $T'$, and
let $\mathcal{M}$ be  any model of $T$ omitting all types of $P$.
Then since $\mathcal{M}$ is a model of $T'$, every element
satisfying $\phi$ must realizes $p$; hence no element of
$\mathcal{M}$ satisfies $\phi$. So $\neg \exists \bar x \phi(\bar
x)$ is in $T'$. Contradiction.

Clearly if $q \in Q$ is supported over $T'$ then $q$ is supported
over $\mathcal{K}(T|P)$. Therefore, by the standard omitting types
theorem, there is a countable model $\mathcal M$ of $T'$ omitting
all types of $P$ and all types of $Q$. Since $T \subseteq T'$,
$\cal M$ is a model of $T$ and thus $\mathcal M$ is in
$\mathcal{K}(T|P)$ and it omits every $q \in Q$. \qed

\bigskip
We note that to derive the standard omitting theorem from the
above theorem, it is sufficient to take $P$ to be the empty set.

Let $L_n(\bar c)=L \cup \{\bar c\}$,  where $\bar c$ is a new
$n$-tuple of constants symbols. For an $n$-generated model $\cal
M$ of $T$, generated by $\bar a$, we let
$$ p_{\bar a}(\mathcal M)=\{
\phi(\bar c)~|~\phi \hbox{ is formula such that } \mathcal M
\models \phi(\bar a)\}, $$
 in the language $L_n(\bar c)$.
We let, as in the introduction,
$$
p_n(y)=\{y \neq  \tau (\bar c)~|~\tau(\bar x)\hbox{ is a term in }
L, ~|\bar x|=n\}.
$$

As noticed in the introduction,  a model $(\mathcal M, \bar a)$ of
$T$, in the language $L_n(\bar c)$, is generated by $\bar a$ if
and only $(\mathcal M, \bar a)$ omits $p_n$.

\begin{defn}
 Let $\mathcal{M}$ be an $n$-generated model of $T$.
 We define inductively $rk_n(\mathcal M) = \alpha$ for an ordinal  $\alpha \geq
 1$ as follows.

\smallskip
 $\bullet$ $rk_n(\mathcal M)=1$ if and only if for every generating $n$-tuple
 $\bar a$ of $\cal M$,  $p_{\bar a}(\mathcal
M)$ is supported over $\mathcal K (T| P_{n,1})$, where
$P_{n,1}=\{p_{n}\}$.

\smallskip
$\bullet$ $rk_n(\mathcal M) = \alpha$ if and only if for every
generating $n$-tuple $\bar a$ of $\cal M$,  $p_{\bar a}(\mathcal
M)$ is supported over $\mathcal K (T| P_{n, \alpha})$, where for
$\alpha \geq 2$,
$$
P_{n, \alpha} =\{p_{0,n}\} \cup \{p_{\bar b}(\mathcal N)
~|~\mathcal N \models T, \bar b \hbox{ generates } \mathcal N,
|\bar b|=n, ~rk_n(\mathcal N) < \alpha\}.
$$

\smallskip
$\bullet$ $rk_n(\mathcal M) = \infty $ if there is no ordinal
$\alpha$ such that $rk_n(\mathcal M) = \alpha $.
\end{defn}

 It should be
noted that we work in the language $L_n(\bar c)$.

We first show that this rank does not depends on $n$.

\begin{prop}\label{pro-invariance}
Let $T$ be a theory in $L$.  If $\cal M$ is an $n$-generated model
of $T$, which is also $m$-generated, then $rk_n(\mathcal
M)=rk_m(\mathcal M)$.
\end{prop}

\proof

Let $n,m \in \mathbb{N}^*$. We prove by induction on $\gamma \geq
1 $, that if $\cal M$ is an $n$-generated of $T$, which is also
$m$-generated, and $rk_n(\mathcal M)=\gamma$,  then $rk_m(\mathcal
M)=\gamma$.

Let $\cal M$ be an $n$-generated model of $T$, which is also
$m$-generated and suppose that $rk_n(\mathcal M)=\gamma$. Then for
any generating $n$-tuple $\bar a$ of $G$, there exists a sentence
$\phi(\bar c)$, with  $|\bar c|=n$, such that for any
$n$-generated model $\cal A$ of $T$ generated by the $n$-tuple
$\bar h$,  if $\mathcal A \models \phi(\bar h)$ and $(\mathcal A,
\bar h) \in \mathcal K(T|P_{n, \gamma})$, then $(\mathcal M, \bar
a) \cong (\mathcal A, \bar h)$.

Let $\bar b$ be an $m$-tuple which generates $\cal M$. Then there
exist an $n$-tuple $\bar t$ of  terms  such that $\bar a=\bar
t(\bar b)$ in $\cal M$, and an $m$-tuple $\bar \tau$ of  terms
 such that
$\bar b=\bar \tau(\bar a)$ in $\cal M$. Let
$$
\psi(\bar y) = \exists \bar x  (\bar x=\bar t(\bar y)\wedge \bar
y=\bar \tau(\bar x) \wedge \phi(\bar x)).
$$

Then $\mathcal M \models \psi(\bar b)$.

Suppose first  $\gamma=1$.  We claim that $\psi(\bar y)$ supports
$p_{\bar b}(\mathcal M)$ over $\mathcal K(T|P_{m,1})$.  Let
$\mathcal A$ be a model of $T$, generated by the $m$-tuple $\bar
d$, such that $\mathcal A \models \psi(\bar d)$. Since $\mathcal A
\models \psi(\bar d)$, there exists an $n$-tuple $\bar h$, which
generates $\cal A$,  such that $\mathcal A \models \phi(\bar h)$.
Then $(\mathcal A, \bar h) \cong (\mathcal M, \bar a)$. Therefore
$(\mathcal A, \bar d) \cong (\mathcal M, \bar b)$, as $\mathcal A
\models \bar d=\bar \tau(\bar h)$ and $\mathcal A \models \bar
b=\bar \tau(\bar a)$. Hence $rk_m(\mathcal M)=1$, and this ends
the proof of our claim and the proof in the case $\gamma=1$.

Now suppose that $rk_n(\mathcal M)=\gamma \geq 2$ and that for any
ordinal $1 \leq \delta < \gamma$ and for any $n$-generated model
$\mathcal N$ of $T$, if $\mathcal N$ is $m$-generated and  if
$rk_n(\mathcal N)=\delta$ then  $rk_m(\mathcal N)=\delta$.

We claim that $\psi(\bar y)$ supports $p_{\bar b}(\mathcal M)$
over $\mathcal K(T|P_{m, \gamma })$.  Let $\mathcal A$ be a model
of $T$, generated by the $m$-tuple $\bar d$, such that $\mathcal A
\models \psi(\bar d)$ and $(\mathcal A, \bar d) \in \mathcal
K(T|P_{m, \gamma })$, and let us prove that $(\mathcal M, \bar b)
\cong (\mathcal A, \bar d)$.

As before,  there exists an $n$-tuple $\bar h$, which generates
$\cal A$,  such that $\mathcal A \models \phi(\bar h)$. Let us
prove that $(\mathcal A, \bar h) \in \mathcal K(T|P_{n, \gamma})$.
Suppose towards a contradiction, that $(\mathcal A, \bar h) \not
\in \mathcal K(T|P_{n, \gamma})$ and thus $(\mathcal A, \bar h)$
realizes $p_{\bar f}(\mathcal N)$ for some $\mathcal N \models T$
and $rk_n(\mathcal N)=\delta < \gamma$. Then $(\mathcal A, \bar h)
\cong (\mathcal N, \bar f)$ and thus $rk_n(\mathcal A)=
rk_n(\mathcal N)=\delta< \gamma$. Therefore, by induction,
$rk_m(\mathcal A)=rk_n(\mathcal A)< \gamma$; as $\mathcal A$ is
$m$-generated. A contradiction as $(\mathcal A, \bar d) \in
\mathcal K(T|P_{m, \gamma })$. Therefore $(\mathcal A, \bar h) \in
\mathcal K(T|P_{n, \gamma})$.

Since  $(\mathcal A, \bar h)  \in \mathcal K(T|P_{n, \gamma})$ and
$\mathcal A \models \psi(\bar d)$,  $(\mathcal A, \bar h) \cong
(\mathcal M, \bar a)$. Therefore $(\mathcal A, \bar d) \cong
(\mathcal M, \bar b)$, as $\mathcal A \models \bar d=\bar
\tau(\bar h)$ and $\mathcal A \models \bar b=\bar \tau(\bar a)$.
Thus $\psi(\bar y)$ supports $p_{\bar b}(\mathcal M)$ over
$\mathcal K(T|P_{m, \gamma})$ as claimed. Hence $rk_m(\mathcal
M)=\gamma$. This ends the proof of the induction.

Now if $\cal M$ is an $n$-generated model of $T$, which is also
$m$-generated,  with $rk_n(\mathcal M) = \infty$, then by the
precedent  result, we deduce $rk_m(\mathcal M) = \infty$. \qed

\smallskip
Proposition \ref{pro-invariance} allows us to define the rank of
$\mathcal M$ in a way that does not depends on the length of
generating tuple of $\mathcal M$. So we let $rk(\mathcal M)$ to be
$rk_n(\mathcal M)$ if $\mathcal M$ is $n$-generated.  When
$\alpha_n(T) \leq \aleph_0$, we have a good characterization.

\begin{thm}\label{rank-theo}
If $T$ is an $n$-consistent theory satisfying $\alpha_n(T) \leq
\aleph_0$, then every $n$-generated model of $T$ has an ordinal
rank.
\end{thm}

\begin{lem} \label{lem-rank} Let $\cal M$ be an $n$-generated model of $T$, generated by
the $n$-tuple $\bar a$, such that $p_{\bar a}(\mathcal M)$ is
supported over  $\mathcal K(T| P_{n, \gamma})$. Then for every
generating $n$-tuple $\bar b$ of $\mathcal M$ the type $p_{\bar
b}(\mathcal M)$ is also supported over $\mathcal K(T| P_{n,
\gamma})$.
\end{lem}

\proof There is a finite $n$-tuple of terms $\bar \tau(\bar x)$
such that $\mathcal{M} \models \bar b = \bar \tau (\bar a)$. Let
$\phi(\bar c)$ be a sentence supporting $p_{\bar a}(\mathcal M)$
over  $\mathcal K(T| P_{n, \gamma})$, and   Let $ \psi(\overline
c)= \exists \bar  x (\bar c=\bar \tau (\bar x))\wedge \phi(\bar
x)). $ We claim that $\psi(\overline c)$ supports $p_{\bar
b}(\mathcal M)$ over $\mathcal K(T| P_{n, \gamma})$. Clearly
$(\mathcal M, \bar a) \in \mathcal K(T|P_{n, \gamma})$ and
$\mathcal M \models \psi(\bar b)$.

Let $(\mathcal{N}, \bar h) \in \mathcal K (T| P_{n, \gamma})$ such
that  $\mathcal N \models \psi(\overline h)$. Then, there exists
$\bar d$ such that $ \mathcal{N} \models \phi(\bar d)$ and
$\mathcal{N}\models (\bar h=\bar \tau(\bar d))$. Therefore,
$\mathcal N$ is generated by $\bar d$. So $(\mathcal{N}, \bar d)$
omits $p_n$, and clearly $(\mathcal {N}, \bar d)$ is in $\mathcal
K (T| P_{n, \gamma})$. Therefore, as $\mathcal{N} \models
\phi(\bar d)$, $(\mathcal M, \bar a) \cong (\mathcal N, \bar d)$.
Hence, $(\mathcal M, \bar b) \cong (\mathcal N, \bar h)$, as
$\mathcal M \models \bar \tau(\bar a)=\bar b $ and $\mathcal N
\models \bar \tau(\bar d)=\bar h $. Thus $(\mathcal N, \bar h)$
realizes $p_{\bar b}(\mathcal M)$. \qed

\smallskip
\noindent \textbf{Proof of Theorem \ref{rank-theo}}.

Let $\gamma$ be the least ordinal such that if $\cal M$ is an
$n$-generated model of $T$ of ordinal rank, then $rk(\mathcal M) <
\gamma$. Suppose towards a contradiction that there exists an
$n$-generated model $\cal N$ of $T$ such that $rk(\mathcal N)=
\infty$. We are going to  prove that there exists an $n$-generated
model $\cal M$ of $T$ of rank $\gamma$ and thus we get a
contradiction.

Since $\alpha_n(T) \leq \aleph_0$,  $P_{n, \gamma}$ is countable.
Notice  that every $n$-generated model of $T$ is either in
$\mathcal K(T| P_{n, \gamma})$ or it has an  ordinal rank.

By our supposition above, for every generating $n$-tuple $\bar a$
of $\cal N$, the model $(\mathcal N, \bar a)$ omits every $p$ in
$P_{n, \gamma}$. Hence $\mathcal K(T|P_{n, \gamma})$ is not empty.
Now we prove the following claim.

\smallskip
\noindent\emph{Claim}. \emph{There exists an $n$-generated model
$\cal M$ of $T$, generated by an $n$-tuple $\bar a$, such that
$p_{\bar a}(\mathcal M)$ is supported over  $\mathcal K(T| P)$.}

\smallskip
\noindent \emph{Proof.} Since $\alpha_n(T) \leq \aleph_0$, there
exists at most $\aleph_0$ nonisomorphic $n$-generated models of
$T$ which have an infinite rank. Let $((\mathcal M_i, \bar a_i),
~\bar a_i \hbox{ generates } \mathcal M_i : i \in \beta \leq
\aleph_0)$ be the list of nonisomorphic $n$-generated models of
$T$, such that   $rk(\mathcal M_i)=\infty$.

Suppose that for every $i \in \beta$, $p_{\bar a_i}(\mathcal M_i)$
is unsupported over $\mathcal K(T| P_{n, \gamma})$. Then by
Theorem \ref{omi-type}, there exists a model $(\mathcal N, \bar
b)$ of $T$ in $\mathcal K(T| P_{n, \gamma})$ which omits  $p_{\bar
a_i}(\mathcal M_i)$ for every $i \in \beta$.

Therefore $\mathcal N \not \cong \mathcal M_i $ for every $i \in
\beta$, and since  $(\mathcal N, \bar b)$ omits every $p \in P_{n,
\gamma}$ we have $rk(\mathcal N)= \infty$. A contradiction.

Hence there exists $\ell \in \beta $ such that $p_{\bar
a_\ell}(\mathcal M_\ell)$ is supported over  $\mathcal K(T| P_{n,
\gamma})$. \qed

\smallskip
By Lemma \ref{lem-rank} and by the Claim above,  there exists an
$n$-generated model $\mathcal M$ of $T$, such that $rk(\mathcal
M)=\gamma$. A final contradiction. \qed

\bigskip

\noindent \textbf{Remark}.  It should be remarked that $\mathcal M
\in \mathcal K(T|P_{n, \gamma})$ if and only if $rk(\mathcal M)
\geq \gamma$. This property will be used freely without any
reference to it.

\noindent \textbf{Examples.}

(1). Let $\Gamma$ be the universal theory of torsion-free abelian
groups. Then every finitely generated model of $\Gamma$ is free
abelian of finite rank.  In fact, it is well known that for any $n
\in \mathbb N^*$,  every finitely generated  group which satisfies
$Th(\mathbb Z^n$ is isomorphic to $\mathbb Z^n$. Therefore for
every free abelian group $G$ of finite rank $n$, we have
$rk(G)=1$, relatively to $Th(\mathbb Z^n)$. We will see that the
rank of nontrivial free abelian groups coincide with another rank
$Rk$ defined for universal theories, relatively to $\Gamma$.
(section 5).

(2). If $G$  is a finitely generated abelian group  and if we let
$\Gamma$ to be the complete theory of $G$, then every finitely
generated model of $\Gamma$ is isomorphic to $G$ and thus
$\alpha(\Gamma)=rk(G)=1$. According to \cite{Nies}, a finitely
generated group $G$ is said \emph{quasi-finitely axiomatisable}
(abbreviated QFA), if there exists a sentence $\phi$ satisfied by
$G$ such that any finitely generated group satisfying $\phi$ is
isomorphic to $G$. A. Nies \cite{Nies} proves that the free
nilpotent group of class 2 with 2 generators is QFA. F. Oger and
G. Sabbagh \cite{oger-sabbagh} generalize this result by showing
that any finitely generated free nilpotent group of class $\geq 2$
is QFA. Moreover they prove that a finitely generated nilpotent
group is QFA if and only if it is prime model of its theory. Thus
the complete theory $\Gamma$ of a finitely generated free
nilpotent group of class $\geq 2$ satisfies $\alpha(\Gamma)=1$.

It is useful to notice  that in general $\alpha_n(T) \leq
\aleph_0$ does not implies $\alpha_m(T) \leq \aleph_0$, for $m
\geq n$. For instance the group theory $T_{gp}$ has at most
$\aleph_0$ nonisomorphic $1$-generated groups, which are the
cyclic groups. However, it is known that $T_{gp}$ has
$2^{\aleph_0}$ nonisomorphic $2$-generated groups (see
\cite{LyndonSchupp77}).

\section{A property analogue to the  Hopf property of groups}

We give in this section  some properties of $n$-generated models
of a theory $T$ satisfying $\alpha_n(T) \leq \aleph_0$. One of
those properties is  analogue  to the Hopf property of a group.
Recall that a group $G$ is said \emph{Hopfian} or has  the
\emph{Hopf property} if every surjective morphism from $G$ to $G$
is an isomorphism.

\begin{thm} \label{theo-prime} Let $T$ be an $n$-consistent theory in $L$ such that $\alpha_n(T)
\leq \aleph_0$ and let $\cal M$ be an $n$-generated model of $T$.
Then for every generating tuple $\overline a$ of $\mathcal{M}$,
there exists a formula $\phi(\overline x)$ such that for every
finitely generated model $\cal N$ of $T$, generated by $\overline
b$, $(\mathcal{M}, \bar a) \cong (\mathcal {N},\bar b)$ if and
only if $\mathcal{N} \models \phi(\overline b)$ and $rk(\mathcal
N) \geq rk(\mathcal M)$.
\end{thm}

\proof  By Theorem \ref{rank-theo},  $rk(\mathcal M)=\gamma$ for
some ordinal $\gamma \geq 1$. By the definition of the rank, for
every generating $m$-tuple $\bar a$ of $\cal M$, $p_{\bar
a}(\mathcal M)$ is supported over $\mathcal K (T| P_{m, \gamma})$.
Therefore, there exists a formula $\phi(\overline x)$ such that
$\cal M \models \phi(\bar a)$ and for every $m$-generated model
$\cal N$ of $T$, generated by $\overline b$, if $\mathcal N
\models \phi(\bar b)$ and $(\mathcal N, \bar b) \in \mathcal
K(T|P_{m, \gamma})$ then $(\mathcal N, \bar b)$ realizes $p_{\bar
a}(\mathcal M)$. Thus we get $(\mathcal{M}, \bar a) \cong
(\mathcal {N},\bar b)$ if and only if $\mathcal{N} \models
\phi(\overline b)$. Since $(\mathcal N, \bar b) \in \mathcal
K(T|P_{m, \gamma})$ if and only $rk(\mathcal N) \geq \gamma$, we
get the desired result. \qed

\begin{cor} \label{cor-prime} Let $T$ be an $n$-consistent theory in $L$ such that $\alpha_n(T)
\leq \aleph_0$. Then there exists  an $n$-generated model
$\mathcal{M}$ of  $T$, such that for every generating tuple
$\overline a$ of $\mathcal{M}$, there exists a formula
$\phi(\overline x)$ such that for every finitely generated model
$\cal N$ of $T$, generated by  $\overline b$,  $(\mathcal{M}, \bar
a) \cong (\mathcal {N},\bar b)$ if and only if $\mathcal{N}
\models \phi(\overline b)$.
\end{cor}

\proof By Theorem \ref{theo-prime},  every $n$-generated model of
$T$ of ordinal rank 1, satisfies the conclusions of the corollary.
Let $\gamma$ be the least ordinal such that  there exists an
$n$-generated model of $T$ of ordinal rank $\gamma$. Then, by
definition of the rank, we have $\gamma=1$. \qed

\bigskip
The following theorem translates an \emph{"internal"} property of
all $n$-generated models of a theory $T$ having at most $\aleph_0$
nonisomorphic $n$-generated models.
\begin{thm} \label{theo-formule-prime} Let $T$ be an $n$-consistent theory in $L$ such that $\alpha_n(T)
\leq \aleph_0$, and let $\mathcal{M}$ be an $n$-generated model of
$T$. Then for every generating tuple $\overline a$ of
$\mathcal{M}$ there exists a formula $\phi(\overline x)$ such that
for every generating tuple $\overline b$ of $\mathcal{M}$ we have
$(\mathcal{M}, \bar a) \cong (\mathcal {M},\bar b)$ if and only if
$\mathcal{M} \models \phi(\overline b)$.
\end{thm}

\proof A consequence of Theorem \ref{theo-prime}. \qed

\smallskip
A consequence of the above property of $\cal M$ can be expressed
as follows. If $\bar a$ is an $n$-generating tuple of $\cal M$ and
$f : \mathcal M \rightarrow \mathcal M$ is a surjective morphism
such that $\mathcal M \models \phi(f(\bar a))$, then $f$ is an
isomorphism. When $G$ is a finitely generated Hopfian group, the
formula $\phi(\bar x)$ can be taken to be $\bar x= \bar x$.

We define a finitely generated model $\cal M$ to be
\emph{weak-Hopfian} if for any generating tuple $\bar a$ of $\cal
M$, there exists a formula $\phi(\bar x)$ such that for any
surjective morphism $f : \mathcal M \rightarrow \mathcal M$ if
$\mathcal M \models \phi(f(\bar a))$, then $f$ is an isomorphism.
The next corollary is therefore a consequence of Theorem
\ref{theo-formule-prime}.

\begin{cor} Let $T$ be an $n$-consistent theory in $L$ such that $\alpha(T)
\leq \aleph_0$. Then every finitely generated model of $T$ is
weak-Hopfian. \qed
\end{cor}

A natural question arises in this context. When a finitely
generated model can be prime of its complete theory ? The
following theorem gives a necessary and sufficient condition.

\begin{thm} \label{theorem-prime-model}Let $T$ be  a complete theory
such that $\alpha_n(T) \leq \aleph_0$. An $n$-generated model
$\mathcal{M}$  of $T$ is prime if and only if there exists a
formula $\theta(\bar x)$, satisfied by some tuple in
$\mathcal{M}$, such that if $\mathcal{M}  \models \theta(\bar a)$
then $\bar a$ generates $\mathcal{M}$.
\end{thm}
We use the following classical result.

\begin{prop} \label{prop}{\rm \cite{Hodges(book)93}} Let $\mathcal{M}$ be a countable model. Then $\mathcal{M}$ is a prime model of
its theory iff for every $m \in \mathbb{N}^*$, each orbit under
the action of $Aut(\mathcal{M})$ on $\mathcal{M}^m$ is first-order
definable without parameters.
\end{prop}

\noindent \textbf{Proof of Theorem \ref{theorem-prime-model}} $\,$

 Suppose that $\mathcal{M}$ is a prime model
of its theory and let $\bar a$ generates $\mathcal{M}$. Then there
is some orbit $\mathcal{O}_n$ containing $\bar a$. By Proposition
\ref{prop}, $\mathcal{O}_n$ is definable by a first order formula
$\theta(\bar x)$. Now if $\mathcal{M} \models \theta(\bar b)$,
then there is an automorphism $f$ such that $f(\bar a)=\bar b$,
and therefore  $\bar b$  generates $\mathcal{M}$.

Suppose now  that there exists a formula $\theta(\bar x)$
consistent in $\mathcal{M}$  such that if $\mathcal{M}  \models
\theta(\bar a)$ then $\bar a$  generates $\mathcal{M}$. Let $\bar
a$ in $\mathcal{M}$ such that $\mathcal{M} \models \theta (\bar
a)$. Then by Theorem \ref{theo-formule-prime}, there exists a
sentence $\phi(\bar c)$ in $L(\bar c)$ such that $\mathcal{M}
\models \phi(\bar a)$, and if $\bar b$ generates $\mathcal{M}$
such that $\mathcal{M} \models \phi(\bar b)$, then  the function
defined by $f(\bar a)=\bar b$  extends to an automorphism.

\smallskip
Let $\mathcal{O}_m$ be an orbit and let $\bar t$ be an $m$-tuple
of terms such that $\bar t_1(\bar a) \in \mathcal {O}_m$. Let us
show that $\mathcal{O}_m$ is defined by the  formula
$$
\psi(\bar y)=\exists \bar z (\phi(\bar z)\wedge\theta(\bar
z)\wedge  \bar y=\bar t(\bar  z)).
$$

Let $\bar b \in \mathcal{O}_m$. Then there is an automorphism $f$
such that $f(\bar t(\bar a))=\bar b$.  Therefore $\bar t(f(\bar
a))=\bar b$, and $\mathcal{M} \models \phi(f(\bar a))\wedge \theta
(f(\bar a))$. Thus $\mathcal{M}  \models \psi(\bar b)$.

Now let $\bar b \in \mathcal{M}$ such that $\mathcal{M} \models
\psi(\bar b)$. Then there is tuple $\bar d$ in $\mathcal{M}$ such
that $\mathcal{M}  \models \phi(\bar d)\wedge\theta(\bar d)\wedge
\bar b=\bar t(\bar d)$. Hence $\bar d$ generates $\mathcal{M}$ and
since $\mathcal{M} \models \phi(\bar d)$ there is an automorphism
$f$ such that $f(\bar a)=\bar d$.  Therefore $f(\bar t(\bar
a))=\bar t(\bar d)=\bar b$. Thus $\bar b \in \mathcal{O}_m$. \qed

\begin{defn} A theory $T$ is said to be
\textbf{$n$-categorical} if $\alpha_n(T)=1$. \rm
\end{defn}

\noindent \textbf{Examples}$ \, \ $

(1)  Let $F_2$ be the free non-abelian group on two generators.
Then  $Th_{\forall \exists }(F_2)$ is $2$-categorical. In fact,
$Th_{\forall}(F_2) \cup \{\exists x\exists y([x,y]\neq 1)\}$ is
$2$-categorical. Indeed, if $A$ is a model of $Th_{\forall}(F_2)$
generated by $\{a,b\}$ and $[a,b] \neq 1$, then it is well-known
that $\{a,b\}$ generates a free group with basis $\{a,b\}$.

(2) For every $n$, $Th(\mathbb {Z}^n)$ is $m$-categorical for
every $m$.

\begin{cor} Let $T$ be  a complete $n$-consistent theory satisfying  $\alpha_n(T) \leq
\aleph_0$. Then the following properties are equivalents

(1)  $\alpha(T)=1$ and the unique finitely generated model of $T$
is prime.

(2) There exists a formula $\theta(\bar x)$, consistent with $T$
such that: for every  finitely generated model $\mathcal{M}$ of
$T$ and for every $\bar b$ in $\mathcal{M}$ if $\mathcal{M}
\models \theta(\bar b)$ then $\bar b$  generates $\mathcal{M}$.

\qed
\end{cor}

\section{The special case of universal theories}

In this section we define another rank specific to universal
theories. In fact this rank can also be defined for $\forall
\exists$-theories and  therefore we work in this context. Before
proceeding, we need some adaptation of Theorem \ref{omi-type} and
Theorem \ref{omi-type2} to our context.

Let $\mathcal{K}$ be a class of models. For an \emph{universal}
type $q$, we say that $q$ is \emph{existentially supported over
$\mathcal{K}$} if there exists an \emph{existential} formula
$\phi(\bar x)$ which supports $q$ over $\cal K$.  We say that $q$
is \emph{existentially unsupported} over $\cal K$ if it is not
existentially supported over $\cal K$.

\smallskip
Let $T$ be a $\forall \exists$-theory in $L$ and $P$ a set of
universal types. The next theorem is a an adaptation of Theorem
\ref{omi-type}. The proof is very similar to the proof of that
theorem.

\begin{thm}\label{omi-type2}
Let $P$ and $Q$ be a countable sets of universal types. If
$\mathcal{K}(T|P)$ is not empty and each $q \in Q$ is
existentially unsupported over $\mathcal{K}(T|P)$, then there
exists a countable model in $\mathcal{K}(T|P)$ which omits every
$q \in Q$.

\end{thm}

\proof

Let $T'$ be the set of all $\forall \exists$-sentences of $L$
which are true in every model in $\mathcal{K}(T|P)$. We claim that
over $T'$, all types of $P$ and all types of $Q$ are with no
existential support; i.e they are existentially unsupported over
the class of models of $T'$.

Let $p \in P$ and suppose that $\phi$ is existential and supports
$p$ over $T'$, and let $\mathcal{M}$ be  any model of $T$ omitting
all types of $P$. Then since $\mathcal{M}$ is a model of $T'$,
every element satisfying $\phi$ must realizes $p$; hence no
element of $\mathcal{M}$ satisfies $\phi$. So $\neg \exists \bar x
\phi(\bar x)$, which is an universal sentence,  is in $T'$.
Contradiction.

Clearly if $q \in Q$ is existential supported over $T'$ then $q$
is existentially supported over $\mathcal{K}(T|P)$. Therefore, by
the  omitting types theorem for $\forall \exists$-theories (see
for instance \cite{Cahng-Keisler} or \cite{Hodges(book)93}), there
is a countable existentially closed model $\mathcal M$ of $T'$
omitting all types of $P$ and all types of $Q$. Since $T \subseteq
T'$, $\cal M$ is a countable model of $T$ and thus $\mathcal M$ is
in $\mathcal{K}(T|P)$ and it omits every $q \in Q$. \qed

\bigskip

If $\cal M$ is a model of $T$, generated by the $n$-tuple $\bar
a$, we let
$$
p_{\bar a, \forall}(\mathcal M)=\{\psi(\bar c)~|~\phi(\bar x)
\hbox{ is universal }, \mathcal M \models \phi(\bar a)\}.
$$

\begin{defn}
 Let $\mathcal{M}$ be an $n$-generated model of $T$.
 We define inductively $Rk_n(\mathcal M) = \alpha$ for an ordinal $\alpha \geq
 1$ as follows.

\smallskip
 $\bullet$ $Rk_n(\mathcal M)=1$ if and only if for every generating $n$-tuple
 $\bar a$ of $\cal M$,  $p_{\bar a}(\mathcal
M)$ is supported over $\mathcal K (T| Q_{n,1})$, where
$Q_{n,1}=\{p_n\}$.

\smallskip
$\bullet$ $Rk_n(\mathcal M) = \alpha$ if and only if for every
generating $n$-tuple $\bar a$ of $\cal M$,  $p_{\bar a}(\mathcal
M)$ is supported over $\mathcal K (T| Q_{n, \alpha})$, where for $
\alpha \geq 2 $
$$
Q_{n, \alpha} =\{p_{0,n}\} \cup \{p_{\bar b, \forall}(\mathcal N)
~|~\mathcal N \models T, \bar b \hbox{ generates } \mathcal N,
|\bar b|=n, ~Rk_n(\mathcal N) < \alpha\}.
$$

\smallskip
$\bullet$ $Rk_n(\mathcal M) = \infty $ if there is no ordinal
$\alpha$ such that $Rk_n(\mathcal M) = \alpha $.
\end{defn}

As in section 3, this rank does not depends on $n$ and we denote
it $Rk$.

Using Theorem \ref{omi-type2}, the proof of the following theorem
is just an adaptation of the proof of Theorem \ref{rank-theo}.
\begin{thm}\label{rank-theo2}
If $T$ is an $n$-consistent $\forall \exists$-theory satisfying
$\alpha_n(T) \leq \aleph_0$, then every $n$-generated model of $T$
has an ordinal rank. \qed
\end{thm}

A natural problem in this context is to explicit the relation
between the two ranks $rk$ and $Rk$.

\begin{prop} If $T$ is an $n$-consistent $\forall \exists$-theory, then for every $n$-generated model
$\mathcal M$ of $T$, $rk(\mathcal M) \leq Rk(\mathcal M)$.

\end{prop}

\proof $\;$

We prove by induction on  $ \gamma \geq 1$, that if $Rk(\mathcal
M)=\gamma$ then $rk(\mathcal M) \leq \gamma$.

For $\gamma=1$ the result follows from the definition of the two
ranks.

Now suppose that for every $n$-generated model $\mathcal N$, if
$Rk(\mathcal N)=\alpha < \gamma $ then $rk(\mathcal M) \leq
\alpha$ and let $\mathcal M$ be an $n$-generated model of $T$ with
$Rk(\mathcal M)= \gamma$.

If $rk(\mathcal M) < \gamma$ we get the searched result. So we
suppose  $rk(\mathcal M) \geq  \gamma$ and we show that
$rk(\mathcal M) =  \gamma$.

Let $\bar a$ be an $n$-generating tuple of $\mathcal M$. Then
there exists a formula $\phi(\bar c)$ which supports $p_{\bar a,
\forall}(\mathcal M)$ over $\mathcal K(T|Q_{n, \gamma})$. We claim
that $\phi(\bar c)$ supports $p_{\bar a}(\mathcal M)$ over
$\mathcal K(T|P_{n, \gamma})$.

Let $\cal N$ be an $n$-generated  model of $T$, generated by $\bar
b$, such that $\mathcal N \models \phi(\bar b)$ and $rk(\mathcal
N) \geq \gamma$. Then  $Rk(\mathcal N) \geq \gamma$; otherwise if
$Rk(\mathcal N) < \gamma$, by induction $rk(\mathcal N) < \gamma$,
a contradiction. Thus $(\cal N, \bar b)$ omits every $p \in Q_{n,
\gamma}$ and hence $(\mathcal N, \bar b) \cong (\mathcal M, \bar
a)$.

We conclude that if $\cal N$ is an $n$-generated  model of $T$,
generated by $\bar b$, such that $\mathcal N \models \phi(\bar b)$
and $rk(\mathcal N) \geq \gamma$, then $(\mathcal N, \bar b) \cong
(\mathcal M, \bar a)$. As $rk(\mathcal M) \geq \gamma$, $\mathcal
K(T|P_{n, \gamma})$ is not empty and therefore $\phi(\bar c)$
supports $p_{\bar a}(\mathcal M)$ over $\mathcal K(T|P_{n,
\gamma})$ as claimed. Thus $rk(\mathcal M)=\gamma$. This ends the
proof of the induction.

If $Rk(\mathcal M)=\infty$ the conclusion is clear. \qed

\smallskip
\noindent \textbf{Remark}.  It is not always the case that
$rk(\mathcal M)=Rk(\mathcal M)$. Let, as in the end of section 3,
$\Gamma$ be the universal theory of torsion-free abelian groups.
Then every finitely generated model of $\Gamma$ is free abelian of
finite rank.  This rank coincide with the above rank; that is, a
nontrivial finitely generated torsion-free abelian group $G$ is of
rank $n$ if and only if $Rk(G)=n$. This shows that, in general,
$rk(\mathcal M)\neq Rk(\mathcal M)$.

\smallskip
The following theorem will be used in the next section.

\begin{thm} \label{prop2}  Let $T$ be an $n$-consistent $\forall \exists$-theory satisfying
$\alpha_n(T) \leq \aleph_0$ and let $\cal M$ be  an $n$-generated
model of $T$, generated by $\bar a$.  Then there exists a
quantifier-free formula $\phi(\bar x)$ such that if $\cal N$ is a
finitely generated  model of $T$, generated by  $\bar b$,  then
$(\mathcal M, \bar a) \cong (\mathcal N, \bar b)$ if and only if
$\mathcal N \models \phi(\bar b)$ and  $Rk(\mathcal N) \geq
Rk(\mathcal M)$.
\end{thm}

\proof

 Let $Rk(\mathcal M)=\gamma$.  By the definition of the
rank, for every generating $n$-tuple $\bar a$ of $\cal M$,
$p_{\bar a, \forall}(\mathcal M)$ is supported over $\mathcal K
(T| Q_{n, \gamma})$. Therefore, there exists an existential
formula $\phi(\overline x)$ such that $\cal M \models \phi(\bar
a)$ and for every $n$-generated model $\cal N$ of $T$, generated
by $\overline b$, if $\mathcal N \models \phi(\bar b)$ and
$(\mathcal N, \bar b) \in \mathcal K(T|Q_{n, \gamma})$ then
$(\mathcal N, \bar b)$ realizes $p_{\bar a, \forall}(\mathcal M)$.

Put $\phi(\bar x) =\exists \bar y \psi(\bar x, \bar y)$, where
$\psi(\bar x)$ is quantifier-free. Since $\mathcal M \models
\psi(\bar a)$, there exists a tuple $\bar d$ such $\mathcal M
\models \psi(\bar a, \bar d)$. Then there exists a tuple of  terms
$\bar t(\bar x)$ such that $\mathcal M \models \bar d=\bar t(\bar
a)$. We let
$$
\xi(\bar x) = \psi(\bar x, \bar t(\bar x)).
$$

 Thus we get $(\mathcal{M}, \bar a) \cong (\mathcal {N},\bar b)$
if and only if $\mathcal{N} \models \xi(\overline a)$. Since
$(\mathcal N, \bar b) \in \mathcal K(T|Q_{n, \gamma})$ if and only
if $Rk(\mathcal N) \geq \gamma$, we get the desired result. \qed

\section{Limit groups of Equationally noetherian groups}

In this section we discuss  properties of limit groups of
equationally noetherian groups, related to their rank.

We begin by recalling some definitions. Let $G$ be a fixed group
and $\bar x=(x_1, \dots, x_n)$. We denote by $G[\bar x]$ the group
${G*F(\bar x)}$ where $F(\bar x)$ is the free group with basis
$\{x_1, \dots,x_n\}$. For an element $s(\bar x) \in G[\bar x]$ and
a tuple  $\bar g =(g_1, \dots, g_n)\in G^n$ we denote by $s(\bar
g)$ the element of $G$ obtained by replacing each $x_i$ by $g_i$
($1 \leq i \leq n$). Let $S$ be a subset of $G[\bar x]$. Then the
set
$$V(S)=\{\bar g \in G^n ~|~s(\bar g)=1 \hbox{ for all }s \in S\} $$
is termed the \emph{algebraic set} over $G$ defined by $S$.  A
group $G$ is called \emph{equationally noetherian} if for every $n
\geq 1$ and every subset $S$ of $G[\bar x]$ there exists a finite
subset $S_0 \subseteq S$ such that ${V(S)=V(S_0)}$.

Let $H$ be a group. For our purpose, we do not need the exact
definition of finitely generated $H$-limit groups, we use an
equivalent definition true
 when $H$ is equationally noetherian
\cite{ould-equa}. A \emph{finitely generated group} $G$ is said
\emph{$H$-limit} if $G$ is a model of the universal theory of $H$.

As noticed in the introduction, if $H$ is an equationally
noetherian group, then the universal theory of $H$ has at most
$\aleph_0$ nonisomorphic finitely generated models. For
completeness we provide a proof of this property.

Throughout this section,  if $G$ is a finitely generated group,
generated by $\bar a$, we let
$$
P(\bar x)=\{ w(\bar  x)~|~w \hbox{ is a word such that } G \models
w(\bar a)=1\}.
$$

\begin{prop} \label{theo-number-finitely} \cite{ould-equa}
Let $H$ be an equationally noetherian group. Then there exist at
most  countably many nonisomorphic finitely generated $H$-limit
groups.
\end{prop}

\proof

Suppose towards a contradiction that the opposite is true. Then
there exists $n \in \mathbb N$ such that there exists at least
$\lambda$ nonisomorphic $n$-generated $H$-limit groups for some
$\lambda > \aleph_0$. Let $(G_i=\<\bar x | P_i(\bar x) \>|i \in
\lambda
> \aleph_0)$ be the list of nonisomorphic $n$-generated
$H$-limit groups. For every $i \in \lambda$ there exists a finite
subset $S_i \subseteq P_i$ such that $ H \models \forall \bar
x(S_i(\bar x)=1 \Rightarrow w(\bar x)=1)$ for every $w \in
P_i(\bar x)$.

Since for every $i \in  \lambda$ the set $S_i$ is finite,  the set
$\{S_i |i \in \lambda\}$ is countable. Therefore the map $f :
\{P_i | i \in \lambda\}\rightarrow \{S_i |i \in \lambda\}$ defined
by $P_i \mapsto S_i$ is not injective and thus there exist $i, j
\in \lambda, i \neq j$ such that $S_i=S_j$.

Since $G_i, G_j$ are models of the universal theory of $H$ we get
$P_i=P_j$, a contradiction. \qed

Therefore,  by Theorem \ref{rank-theo2}, every finitely generated
$H$-limit group $G$ has an ordinal rank $Rk$. We are interested on
the relation between this rank and decomposition of morphisms from
$G$ to another $H$-limit group.

\begin{thm}\label{theo-finite-quotient} Let $H$ be an equationally noetherian
group and $G$ be a nontrivial  finitely generated $H$-limit group.
Then there exists a finite collection of proper epimorphisms of
$H$-limit groups $(f_i : G \rightarrow L_i|~1 \leq i \leq n)$,
such that $Rk(L_i)\geq Rk(G)$ whenever $L_i$ is nontrivial, and
such that for any $H$-limit group $L$, if $Rk(L) \geq Rk(G)$, then
any epimorphism $f :G\rightarrow L$, is either an embedding or
factors through some
 $f_i$.
\end{thm}

\proof

Let $\bar a$ be an $n$-tuple which generates $G$. Set
$\alpha=Rk(G)$. Let $(f_i : G \rightarrow G_i~|~i \in \N)$ be the
list of all
  proper quotients of $G$ which are $H$-limit (Notice that the trivial group is a proper quotient of $G$). Then every $G_i$ is generated
 by $f_i(\bar a)$. As before, we let
 $$
 P(\bar x)=\{
w(\bar  x)~|~w \hbox{ is a word such that } G \models w(\bar
a)=1\},
 $$
 $$
P_i(\bar x)=\{ w(\bar  x)~|~w \hbox{ is a word such that } G_i
\models w(f_i(\bar a))=1\}.
 $$

Since $H$ is equationally noetherian,  there exist a finite
subsets $S(\bar x) \subseteq P(\bar x)$, $S_i(\bar x) \subseteq
P_i(\bar x)$ such that
$$
H \models \forall \bar x(S(\bar x)=1 \Rightarrow w(\bar x)=1),
\hbox{ for any }w \in P(\bar x),
$$
$$
H \models \forall \bar x(S_i(\bar x)=1 \Rightarrow w(\bar x)=1),
\hbox{ for any }w \in P_i(\bar x).
$$

By Theorem \ref{prop2}, for any $G_i$,  there exists a
quantifier-free formula $\phi_{i}(\bar  x)$ such that if $A$ is an
$H$-limit group, if  $Rk(A) \geq Rk(G_i)$ and if $\bar d$ is an
$n$-generating tuple  of $A$ such that $A \models \phi_{i}(\bar
d)$, then $(G_i, f_i(\bar a)) \cong (A, \bar d)$.

Similarly, for $G$ and $\bar a$, there exists a quantifier-free
formula $\phi(\bar x)$ such that if $A$ is an $H$-limit group such
that $Rk(A) \geq Rk(G)$ and if $\bar d$ is an $n$-generating tuple
of $A$ such that $A \models \phi(\bar d)$, then $(G, \bar a) \cong
(A, \bar d)$.

\smallskip
Let
$$
\Gamma(\bar c)=Th_{\forall}(H) \cup \{\neg \phi_{i }(\bar c) |
~Rk(G_i) <Rk(G)\} \cup \{\neg \phi(\bar c)\} \cup \{S(\bar c)=1\}.
$$

If $\Gamma(\bar c)$ is not consistent, then
$$
Th_{\forall}(H) \cup \{\neg \phi_{i }(\bar c) | ~Rk(G_i) <Rk(G)\}
\cup \{S(\bar c)=1\} \vdash \phi(\bar c),
$$
and thus, if  $f : G \to  L$ is an $H$-limit quotient with $Rk(L)
\geq \gamma$, then $L \models \phi(f(\bar a))$, and thus  $f$ is
an embedding. Therefore by taking $(f : G \to 1)$ to be our
sequence, we get the result.

So we suppose that $\Gamma (\bar c)$ is consistent. We claim that
$$
\Gamma(\bar c) \vdash \bigvee_{i \in I} S_i(\bar c)=1, \hbox{
where } I=\{ i \in \omega | Rk(G_i) \geq Rk(G)\}. \leqno (1)
$$

Let $(\mathcal M, \bar d)$ be a model of $\Gamma(\bar c)$ in the
language $L(\bar c)$. Let $A$ be the subgroup of $\mathcal M$,
generated by $\bar d$. Then $A$ is an $H$-limit group and since $A
\models S(\bar d)=1$, there exists a morphism $f$ from $G$ to $A$
whcih sends $\bar a$ to $\bar d$. Also since $A \models \neg
\phi(\bar d)$, $A$ is a proper quotient of $G$. Furtheremore,
$Rk(A) \geq Rk(G)$; because if $Rk(A) < Rk(G)$ then $A=G_i$ and $A
\models \phi_i(\bar d)$ for some $i$, a contradiction.  This ends
the proof of our claim.

By compactness and $(1)$,  we get
$$
\Gamma(\bar c) \vdash S_{i_1}(\bar c)=1 \vee \cdots \vee
S_{i_m}(\bar c)=1.
$$

Let $L_j=G_{i_j}$ for $1 \leq j \leq m$ and $(f_i : G \rightarrow
L_i|~1 \leq i \leq m)$ defined obviously. Then this sequence
satisfy the desired conclusion. \qed

\smallskip
The following theorem is a generalization of \cite[Theorem
2.6]{ould-equa}
\begin{thm} \label{thm-gamma} Let $H$ be an equationally noetherian
group and $G$ be a finitely generated $H$-limit group. Then for
any finite subset $X \subseteq G\setminus \{1\}$ and for any
ordinal $1 \leq \gamma \leq Rk(G)$, there exists an epimorphism $f
: G \rightarrow L$ such that $L$ is an $H$-limit group with
$Rk(L)= \gamma$ and $1 \not \in f(X)$.
\end{thm}

Before the proof we need some notions and results from
\cite{ould-equa}.

\begin{defn}
A finitely generated   $H$-limit group $G$ is said
\emph{$H$-determined} if  there exists a finite subset $X
\subseteq G \setminus \{1\}$ such that for any morphism $f : G
\rightarrow L$, where $L$ is an $H$-limit group, if $1 \not \in
f(X)$ then $f$ is an embedding.
\end{defn}

\begin{lem} \label{lem-h-deter} Let $H$ be an equationally noetherian group. A
finitely generated $H$-limit group $G$ is $H$-determined if and
only if $Rk(G)=1$.
\end{lem}

\proof Let $G$ be a finitely generated $H$-limit group which is
also $H$-determined and let us prove that $Rk(G)=1$. Write
$G=\<\bar a|P(\bar a)\>$. Then  there exists a finite subset
$S(\bar x) \subseteq P(\bar x)$ such that $ H \models \forall \bar
x(S(\bar x)=1 \Rightarrow w(\bar x)=1), \hbox{ for any }w \in
P(\bar a). $ Let $X$ given by the words $v_1(\bar x), \cdots,
v_m(\bar x)$. We claim that the formula
$$
\phi(\bar c) \equiv (S(\bar c)=1 \wedge \bigwedge_{1 \leq i \leq
m} v_i(\bar x) \neq 1),
$$
supports $p_{\bar a, \forall}(G)$ over $P_{n,1}$.

Let $L$ be an $H$-limit group, generated by $\bar b$, such that $L
\models \phi(\bar b)$. Then there exists a morphism $f : G \to L$
which sends $\bar a$ to $\bar b$ and $f(v_i(\bar a)) \neq 1$.
Since $G$ is $H$-determined we find that $f$ is an isomorphism and
thus $(G, \bar a) \cong (L, \bar b)$. Therefore $\phi(\bar c)$
supports $p_{\bar a, \forall}(G)$ over $P_{n,1}$ and thus
$Rk(G)=1$.

Now suppose that $Rk(G)=1$ and let us prove that $G$ is
$H$-determined. By theorem \ref{prop2},  there exists a
quantifier-free formula $\phi(\bar x)$ such that if $L$ is an
$H$-limit group, generated by  $\bar b$, then $(G, \bar a) \cong
(L, \bar b)$ if and only if $L \models \phi(\bar b)$ and $Rk(L)
\geq Rk(G)=1$. Clearly, replacing $\phi$ by a primitive
quantifier-free formula, one can assume that $\phi$ is primitive
quantifier-free. Then
$$
\phi(\bar x) \equiv (\bigwedge_{w \in W} w(\bar x)=1 \wedge
\bigwedge_{v \in V} v(\bar x) \neq 1),
$$
where $W, V$ are finite sets of words.

Therefore, if $f : G \to L$ is a morphism, with $L$ is $H$-limit
and $f(v(\bar a))\neq 1$, for any $v \in V$, then $f$ is an
embedding. Thus $G$ is $H$-determined as desired. \qed

\smallskip
We will also need the following.

\begin{thm} \label{theo-qutient-finitely-determined} \cite{ould-equa} Let $H$ be an equationally noetherian group and $G$  a
nontrivial finitely generated $H$-limit group. Then for any finite
subset $X \subseteq G \setminus \{1\}$ there exists an epimorphism
$f : G \rightarrow L$ where $L$ is an $H$-determined group such
that $1 \not \in f(X)$. \qed

\end{thm}

\smallskip
\noindent \textbf{Proof of Theorem \ref{thm-gamma}}.

If $\gamma=1$, then the theorem is a consequence of Lemma
\ref{lem-h-deter} and Theorem
\ref{theo-qutient-finitely-determined}. So we suppose that $\gamma
\geq 2$.

Suppose that $G$ is $n$-generated. Let $( (G_i, \bar a_i)~|~ \bar
a_i \hbox{ generates } G_i, |\bar a_i|=n, ~i \in \omega)$ be the
list of all
 $H$-limit groups, up to isomorphism,  such that $1 \leq Rk(G_i)<\gamma$.

By theorem \ref{prop2}, for any $(G_i, \bar a_i)$,  there exists a
quantifier-free formula $\phi_{i}(\bar  x)$ such that if $A$ is an
$H$-limit group, if  $Rk(A) \geq Rk(G_i)$ and if $\bar d$ is an
$n$-generating tuple  of $A$ such that $A \models \phi_{i}(\bar
d)$, then $(G_i, \bar a_i) \cong (A, \bar d)$.

Let
$$
\Gamma(\bar c)=Th_{\forall}(H) \cup \{\neg \phi_{i }(\bar c) | ~i
\in \omega\}.
$$
Then $\Gamma (\bar c)$ is consistent as $(G, \bar a)$ is model of
$\Gamma (\bar c)$, for any generating $n$-tuple $\bar a$ of $G$.
Now we prove the following claim.

\smallskip
\noindent \emph{Claim.} \emph{For any primitive-quantifier-free
formula $\vartheta (\bar x)$ such that $\Gamma(\bar c) \cup
\{(\vartheta (\bar c))\}$ is consistent, there exists a
primitive-quantifier-free formula $\xi (\bar x)$ such that
$\Gamma(\bar c) \cup \{(\vartheta (\bar c)\wedge \xi (\bar c))\}$
is consistent  and for any word  $w (\bar x)$ on the variables
$\bar x=\{x_1, \dots, x_n\}$ and their inverses one has
$$
\Gamma(\bar c) \vdash (\vartheta (\bar c)\wedge \xi (\bar
c)\Rightarrow w (\bar c)=1) \hbox{ or } \Gamma (\bar c) \vdash
(\vartheta (\bar c)\wedge \xi (\bar c)\Rightarrow
 w (\bar c) \neq 1).
$$}

\proof

Let $\vartheta (\bar x)$ be a primitive-quantifier-free formula
such that   $\Gamma(\bar c) \cup \{(\vartheta (\bar c))\}$ is
consistent and suppose towards a contradiction that $\vartheta
(\bar x)$ does not satisfies the conclusions of the claim.  We are
going to construct a tree. By hypothesis  there exists  a word
$\alpha _{1}(\bar x)$ such that $\Gamma (\bar c)\cup \{(\vartheta
(\bar c)\wedge \alpha _{1}(\bar c)=1)\}$ and $\Gamma (\bar c)\cup
\{ (\vartheta (\bar x)\wedge \alpha _{1}(\bar c) \neq 1)\}$ are
consistent (to simplify notation we omit $\bar c$). We can do the
same thing with $\vartheta \wedge \alpha _{1}=1$ and $\vartheta
\wedge \alpha _{1}\neq 1$. Thus we have: 
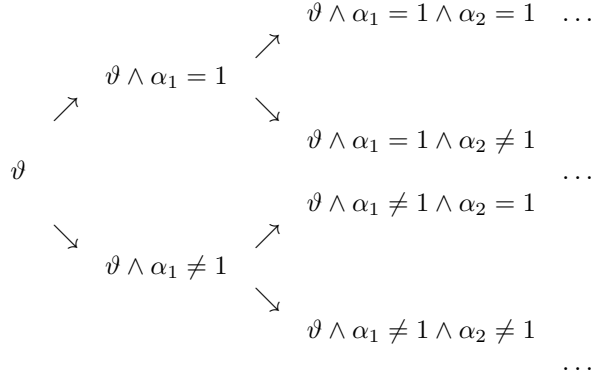
\[
\begin{array}{cccccccc}
  &  &  &  &  & \vartheta \wedge \alpha _{1}=1\wedge \alpha _{2}=1 & \dots & \\
  &  &  &  & \nearrow  &  &  & \\
  &  &  & \vartheta \wedge \alpha _{1}=1 &  &  &  & \\
  &  & \nearrow  &  & \searrow  &  &  & \\
  &  &  &  &  & \vartheta \wedge \alpha _{1}=1\wedge  \alpha _{2}\neq 1 &  & \\
  & \vartheta  &  &  &  &  & \dots & \\
  &  &  &  &  & \vartheta \wedge  \alpha _{1}\neq 1\wedge \alpha _{2}=1 &  & \\
  &  & \searrow  &  & \nearrow  &  &  & \\
  &  &  & \vartheta \wedge  \alpha _{1} \neq 1 &  &  &  & \\
  &  &  &  & \searrow  &  &  & \\
  &  &  &  &  & \vartheta \wedge  \alpha _{1}\neq 1\wedge  \alpha _{2} \neq 1 &  & \\
  &  &  &  &  &  & \dots & \\
  &  &  &  &  &  &  & \end{array}
\]

Therefore, by compactness, every branch in the tree is consistent
with $\Gamma(\bar c)$. Since there exists $2^{\aleph_0}$ branch we
get $2^{\aleph_0}$ nonisomorphic finitely generated models of
$\Gamma(\bar c)$ and thus  Th$_{\forall}(H)$ has $2^{\aleph_0}$
nonisomorphic finitely generated models. A  contradiction. \qed

Write $G=\<\bar a| P(\bar a)\>$, $|\bar a|=n$, and let  $S(\bar a)
\subseteq P(\bar a)$ be a finite set such that
$$
H \models \forall \bar x(S(\bar x)=1 \Rightarrow w(\bar x)=1),
\hbox{ for any }w \in P(\bar x).
$$

Let $X \subseteq G$ be a finite subset, given by  words $v_1(\bar
x), \dots, v_n(\bar x)$ such that $G \models \bigwedge_{1 \leq i
\leq n}v_i(\bar a) \neq 1$. Let $\vartheta(\bar x) = (S(\bar x)=1
\wedge \bigwedge_{1 \leq i \leq n}v_i(\bar x) \neq 1)$.

Since $(G, \bar a) \models \vartheta (\bar a)$, $\Gamma(\bar c)
\cup \{\vartheta (\bar c)\}$ is consistent. Therefore,  by the
claim above, there exists a primitive-quantifier-free formula $\xi
(\bar x)$ such that $\Gamma(\bar c) \cup \{(\vartheta (\bar
c)\wedge \xi (\bar c))\}$ is consistent  and for any word  $w
(\bar x)$ on the variables $\bar x=\{x_1, \dots, x_n\}$ and their
inverses one has
$$
\Gamma(\bar c) \vdash (\vartheta (\bar c)\wedge \xi (\bar
c)\Rightarrow w (\bar c)=1) \hbox{ or } \Gamma (\bar c) \vdash
(\vartheta (\bar c)\wedge \xi (\bar c)\Rightarrow
 w (\bar c) \neq 1). \leqno (*)
$$

Let $(\mathcal M, \bar b)$ a model of $\Gamma(\bar c)  \cup
\{(\vartheta (\bar c)\wedge \xi (\bar c))\}$, and let $L$ be the
subgroup of $\mathcal M$ generated by $\bar b$. We claim that $L$
satisfies the desired property. Clearly,  $L$ is an $H$-limit
group. Since $L \models \vartheta (\bar b)$, there exists an
epimorphism $f : G \to L$ such that $f(\bar a)=\bar b$ and $1 \not
\in f(X)$.

Therefore, it remains to show that $Rk(L)=\gamma$. Since $L
\models \neg \phi_i(\bar b)$ we have $Rk(L) \geq \gamma$. By the
property $(*)$, we see that the formula $\vartheta (\bar x) \wedge
\xi (\bar x)$ supports $p_{\bar b, \forall}(L)$ over $\mathcal K
(Th_{\forall}(H)|Q_{n, \gamma})$ and thus by the definition of the
rank we get $Rk(L)=\gamma$ as desired. \qed

\bibliographystyle{alpha}
\bibliography{biblio}
\end{document}